\documentclass[11pt, leqno]{article}
\usepackage{amssymb}
\usepackage{amsmath}
\usepackage{amsthm}
\usepackage{latexsym}
\usepackage{amsfonts}
\usepackage{graphicx}
\usepackage{graphics}

\theoremstyle{plain}
\theoremstyle{definition}

\newtheorem*{Proof}{Proof}

\newcommand{\el}{\ell}

\newcommand{\ra}{\;\rightarrow\;}

\newcommand{\ga}{\gamma }
\newcommand{\Ga} {{\varGamma}}

\newcommand{\OO} {{\varOmega}}
\newcommand{\De} {{\varDelta}}

\newcommand{\ti}{\tau }

\newcommand{\R}{\mathbb{R}}

\newcommand{\cl}{{\cal{L}}}

\newcommand{\ld}{\ldots}

\newcommand{\sm}{\smallsetminus}

\begin{document}
\title{On Almost Entire Solutions Of The Burgers Equation}
\author{N. D. Alikakos  and D. Gazoulis}
\date{}
\maketitle
{\bf Abstract:} We consider Burger's equation on the whole $x-t$ plane. We require the solution to be classical everywhere, except possibly over a closed set $S$ of potential singularities, which is

(a) a subset of a countable union of ordered graphs of differentiable functions,

(b) has one dimensional Hausdorff measure, $H^1(S)$, equal to zero.

We establish that under these conditions the solution is identically equal to a constant.
\section{Introduction}\label{sec1}
\noindent

In this note we establish a sort of rigidity theorem for solutions of the Burgers equation
\begin{eqnarray}
h_t(x,t)+h(x,t)h_x(x,t)=0   \label{eq1}
\end{eqnarray}
in the plane $\R_x\times\R_t$. We consider functions $h(x,t)$ that solve (\ref{eq1}) classically, pointwise, except perhaps on a closed set $S$ of the $x-t$ plane as in the Abstract, and we show that $h$ must be identically constant. We note that such a statement is false in the half plane $\R_x\times\R^+_t$ because of rarefaction waves. We also note that the conclusion of the theorem is relatively simple to recover for entropy solutions. Indeed if $u(x,t)$ is an $L^\infty(\R_x\times\R_t)$ entropy solution to
\begin{eqnarray}
\left\{\begin{array}{l}
         u_t+\dfrac{1}{2}(u^2)_x=0 \\
         u(x,0)=u_0(x)
       \end{array}\right.  \label{eq2}
\end{eqnarray}
then we have the (well known) estimate
\begin{eqnarray}
\frac{u(x+a,t)-u(x,t)}{a}<\frac{E}{t}  \label{eq3}
\end{eqnarray}
for every $a>0$, $t>0$ with $E$ depending only on $\|u_0\|_{L^\infty}=M$ ([4], Theorem 16-4, or [3], Lemma in 3.4.3). By shifting the origin of time all the way to $t=-\infty$, and by uniqueness in the entropy class, we conclude via (\ref{eq3}) that $x\ra u(x,t)$ is non-\linebreak increasing for every $t$. Thus in particular $u_0$ is a nonincreasing $L^\infty$ function, and if $u_0$ is not identically constant (a.e.) then the solution of (\ref{eq2}) will have a shock. Thus, the hypothesis $H^1(S)=0$ will force $u_0$ to be identically constant, and so also $u$.

There is a similar result for the eikonal equation
\begin{eqnarray}
\bigg(\frac{\partial u}{\partial t}\bigg)^2+\bigg(\frac{\partial u}{\partial x}\bigg)^2=1  \label{eq4}
\end{eqnarray}
by Caffarelli and Crandall \cite{1} which states that if $u$ solves (\ref{eq4}) pointwise on $\R^2\sm\widehat{S}$, with $H^1(\widehat{S})=0$, then necessarily $u$ is either affine or a double ``cone function'', $u(y)=a\pm|y-z|$, $y=(x,t)$, $z=(x_0,t_0)$. The point in \cite{1} again is that $u$ is not assumed a viscosity solution.

The proof of our result is based on a simple and explicit change of variables (see (\ref{eq6}) below) that transforms (\ref{eq2}) into (\ref{eq4}), and actually establishes almost\footnote{Note that for the set $\widehat{S}$ in \cite{1} there is no extra hypothesis besides that $H^1(\widehat{S})=0$.} the equivalence of the two problems in $\R^2$. Our only excuse for writing it down is that it concerns the Burgers equation, which in spite of its simplicity pervades the theory of hyperbolic conservation laws \cite{2}, \cite{3}. We would like to thank the referees for their useful suggestions that improved significantly the presentation of our result.
\section{Theorem}
\noindent

Let $h(x,t)$ be a measurable function on $\R^2$ and suppose that $S$ is closed and on $\R^2\sm S$ the following hold: $h(x,t)$ is continuous, $\dfrac{\partial h}{\partial t}$, $\dfrac{\partial h}{\partial x}$ exist, $x\ra\dfrac{\partial h}{\partial x}(x,t)$ is $L^1_{loc}$ and moreover
\begin{eqnarray}
h_t+hh_x=0, \ \ \text{on} \ \ \R^2\sm S.  \label{eq5}
\end{eqnarray}
If $H^1(S)=0$ and moreover $S\subset\cup_{i\in Z}\Ga_i$, where $\Ga_i:=\{(x,t)\mid t=p_i(x)$, $p_i$ differentiable, $x\in \R\}$, $\ldots<p_{-n}(x)<\cdots<p_{-1}(x)<p_1(x)<p_2(x)<\cdots<p_n(x)<\cdots$ then $h\equiv$\,constant on $\R^2$, and $S=\emptyset$.\bigskip

{\bf Notes}

1) The change of variables $h=c(v)$ converts $v_t+c(v)v_x=0$ into Burgers' equation $h_t+hh_x=0$, hence this more general equation is covered for differentiable $c$ provided that $c'\neq0$. Note that if we write the equation for $v$ in divergence form $v_t+(C(v))_x=0$, where $C'=c$, then the condition $c'\neq0$ corresponds to $C''\neq0$ which is naturally weaker than the usual condition of genuine nonlinearity $C''>0$, since we do not require any orientation of the $x-t$ plane.

2) The change of variables relating (\ref{eq2}) to (\ref{eq4}) is basically
\begin{eqnarray}
u(x,t)=\int^t_0\frac{ds}{\sqrt{h^2(x,s)+1}}+g(x)  \label{eq6}
\end{eqnarray}
where $g(x)=\int\limits^x_0\dfrac{h(u,0)du}{\sqrt{h^2(u,0)+1}}$.\\
Note that the projected characteristics of the corresponding equations coincide,
\[
\begin{array}{cc}
  \dfrac{dx}{d\ti}=h & \dfrac{dx}{d\ti}=u_x=\dfrac{h}{\sqrt{h^2+1}} \\ [2ex]
 \dfrac{dt}{d\ti}=1  & \dfrac{dt}{d\ti}=u_t=\dfrac{1}{\sqrt{h^2+1}}.
\end{array}
\]
The need for differentiating under the integral sign in (\ref{eq6}) for obtaining (\ref{eq4}) forces us to introduce the perhaps unecessary hypothesis that $S$ lies on a set of graphs.

3) The hypotheses on the singular set a priori do not exclude $S$ to be a countable union of Cantor sets arranged on a family of parallel lines in the $x-t$ plane.%\vspace*{0.2cm} \\
%
%\noindent
%
\begin{Proof} For the convenience of the reader we begin by giving the proof in the simple case where $S$ lies on a single differentiable graph contained inside a strip, $S\subset\Ga:=\{(x,t)\mid t=p(x)$, $p$ differentiable, $0<p(x)<1$, $x\in\R\}$.\\
Set $\OO^+=\{(x,t)\in\R^2\mid t\le p(x)\}$, $\OO^-=\{(x,t)\in \R^2\mid t\ge p(x)\}$.\\
Define for $(x,t)\in\OO^+$
\begin{eqnarray}
u^+(x,t)=\int^t_0\frac{ds}{\sqrt{h^2(x,s)+1}}+g^+(x)  \label{eq7}
\end{eqnarray}
where $g^+(x)=\int\limits^x_0\dfrac{h(u,0)du}{\sqrt{h^2(u,0)+1}}$ and for $(x,t)\in\OO^-$
\begin{eqnarray}
u^-(x,t)=\int^t_1\frac{ds}{\sqrt{h^2(x,s)+1}}+g^-(x)  \label{eq8}
\end{eqnarray}
where $g^-(x)=\int\limits^x_0\dfrac{h(u,1)du}{\sqrt{h^2(u,1)+1}}$ \\
We begin with $u^+(x,t)$ for $t\le p(x)$, $(x,t)\in U:=\R^2\sm S$, open.

By our hypothesis
\begin{align}
u^+_x(x,t)&=\int^t_0\frac{-h(x,s)h_x(x,s)}{\big(\sqrt{h^2(x,s)+1}\big)^3}ds+\frac{h(x,0)}{\sqrt{h^2(x,0)+1}} \nonumber\\
&=\int^t_0\frac{h_s(x,s)}{\big(\sqrt{h^2(x,s)+1}\big)^3}ds+\frac{h(x,0)}{\sqrt{h^2(x,0)+1}} \nonumber \\
&=\frac{h(x,t)}{\sqrt{h^2(x,t)+1}}.  \label{eq9}
\end{align}
On the graph we have
\begin{eqnarray}
u^+_x(x,p(x))=\frac{h(x,p(x))}{\sqrt{h^2(x,p(x))+1}}, \ \ (x,p(x))\notin S.  \label{eq10}
\end{eqnarray}
Differentiating in $t$ is straightforward, and holds quite generally,
\begin{eqnarray}
u^+_t(x,t)=\frac{1}{\sqrt{h^2(x,t)+1}}, \ \ u^+_t(x,p(x))=\frac{1}{\sqrt{h^2(x,p(x))+1}}.  \label{eq11}
\end{eqnarray}
Thus from (\ref{eq9}), (\ref{eq11}) we have
\begin{eqnarray}
(u^+_x(x,t))^2+(u^+_t(x,t))^2=1 \ \ \text{in} \ \ \OO^+\sm S.  \label{eq12}
\end{eqnarray}
Analogously we argue for $u^-(x,t)$ and we obtain
\begin{eqnarray}
u^-_x(x,t)=\frac{h(x,t)}{\sqrt{h^2(x,t)+1}} \ \ \text{in} \ \ \OO^-\sm S,  \label{eq13}
\end{eqnarray}
\begin{eqnarray}
u^-_x(x,p(x))=\frac{h(x,p(x))}{\sqrt{h^2(x,p(x))+1}}, \ \ (x,p(x))\notin S,  \label{eq14}
\end{eqnarray}
\begin{eqnarray}
u^-_t(x,t)=\frac{1}{\sqrt{h^2(x,t)+1}}, \ \ u^-_t(x,p(x))=\frac{1}{\sqrt{h^2(x,p(x))+1}}  \label{eq15}
\end{eqnarray}
and so once more
\begin{eqnarray}
(u^-_x(x,t))^2+(u^-_t(x,t))^2=1 \ \ \text{in} \ \ \OO^-\sm S.  \label{eq16}
\end{eqnarray}
Also from (\ref{eq10}), (\ref{eq14}) we obtain
\begin{eqnarray}
u^+_x(x,p(x))=u^-_x(x,p(x)), \ \ u^+_t(x,p(x))=u^-_t(x,p(x)), \ \ (x,p(x))\notin S.  \label{eq17}
\end{eqnarray}
We now set
\begin{eqnarray}
u(x,t)=\left\{\begin{array}{ll}
                u^+(x,t), & (x,t)\in\OO^+ \\
                u^-(x,t)+\De(x), & (x,t)\in\OO^-
              \end{array}\right.  \label{eq18}
\end{eqnarray}
where
\begin{eqnarray}
\De(x):=u^+(x,p(x))-u^-(x,p(x)), \ \ x\in\R.  \label{eq19}
\end{eqnarray}
Note that $\Ga\sm S$ is open in $\Ga$ and so is its projection $\pi_x(\Ga\sm S)=\bigcup\limits^\infty_{i=1}(a_i,b_i)=:O$, and for $x\in O$
\begin{eqnarray}
\frac{d\De(x)}{dx}=u^+_x(x,p(x))+u^+_t(x,p(x))p'(x)-(u^-_x(x,p(x))+u^-_t(x,p(x))p'(x))=0 \label{eq20}\hspace*{-1.5cm}
\end{eqnarray}
(by (\ref{eq17})). Therefore, by the continuity of $h$ and $p$, $u(x,t)$ is differentiable on $\R^2\sm S$, and by (\ref{eq12}), (\ref{eq16}), (\ref{eq18}) and (\ref{eq20})
\begin{eqnarray}
(u_x(x,t))^2+(u_t(x,t))^2=1 \ \ \text{on} \ \ \R^2\sm S.  \label{eq21}
\end{eqnarray}
Hence, by the result in \cite{1} $u$ is of the form
\begin{eqnarray}
u(x,t)=ax+bt+\ga \ \ (a^2+b^2=1),  \label{eq22}
\end{eqnarray}
or
\begin{eqnarray}
u(x,t)=c\pm\sqrt{(x-x_0)^2+(t-t_0)^2}.  \label{eq23}
\end{eqnarray}
In the first case $u_t=b$ and so $h(x,t)\equiv$\,constant.\\
On the other hand (\ref{eq23}) gives
\begin{align}
u_t(x,t)&=\pm\frac{t-t_0}{\sqrt{(x-x_0)^2+(t-t_0)^2}}  \nonumber \\
&\Rightarrow h(x,t)=\frac{x-x_0}{t-t_0}  \label{eq24}
\end{align}
which is singular on $\{t=t_0\}$, and thus is excluded by the hypothesis $H^1(S)=0$. Therefore $h(x,t)\equiv$\,constant is the only option.\vspace*{0.2cm} \\
{\bf Note:} $\De(x)$ is continuous for $x\in\R;\cl(\pi_x(S))=0$.

In the general case we indicate the necessary modifications. Suppose $p_\el(x)<p_{\el+1}(x)$, $a_\el(x)\in C^1$, $p_\el(x)<a_\el(x)<p_{\el+1}(x)$, $\el=1,2,\ld$, $\el=-2,-3,\ld$ (and $p_{-1}(x)<a_0(x)<p_1(x)$) where we have inserted the $C^1$ graphs $a_\el(x)$ that will play the role of the horizontal lines $t=0$ and $t=1$ in the simple case treated above.

Let
\begin{eqnarray}
\OO^+_1=\{p_{-1}(x)\le t\le p_1(x)\}, \ \ \OO^-_1=\{p_1(x)\le t\le a_1(x)\}
\label{eq25}
\end{eqnarray}
\begin{eqnarray}
u^+_1(x,t):=\int^t_{a_0(x)}\frac{ds}{\sqrt{h^2(x,s)+1}}+g^+_1(x), \; g^+_1(x)=\int^x_0\frac{h(s,a_0(s))+a'_0(s)}{\sqrt{h^2(s,a_0(s))+1}}ds, \; \text{on} \; \OO^+_1 \hspace*{-1.7cm} \label{eq26}
\end{eqnarray}
\begin{eqnarray}
u^-_1(x,t):=\int^t_{a_1(x)}\frac{ds}{\sqrt{h^2(x,s)+1}}+g^-_1(x),\;g^-_1(x)=\int^x_0
\frac{h(s,a_1(s))+a'_1(s)}{\sqrt{h^2(s,a_1(s))+1}}ds,\;\text{on}\;\OO^-_1\hspace*{-2cm}   \label{eq27}
\end{eqnarray}
$\De_1(x):=u^+_1(x,p_1(x))-u^-_1(x,p_1(x))$
\begin{eqnarray}
u_1(x,t)=\left\{\begin{array}{lll}
                  u^+_1(x,t), & \text{on} & \OO^+_1 \\
                  u^-_1(x,t)+\De_1(x), & \text{on} & \OO^-_1.
                \end{array}\right.  \label{eq28}
\end{eqnarray}
For $i=2,3,\ld$ set
\begin{eqnarray}
\OO^+_i=\{a_{i-1}(x)\le t\le p_i(x)\}, \ \ \OO^-_i=\{p_i(x)\le t\le a_i(x)\},  \label{eq29}
\end{eqnarray}
\begin{eqnarray}
u^+_i(x,t):=u^-_{i-1}(x,t)+\De_{i-1}(x), \ \ \text{on} \ \ \OO^+_i,  \label{eq30}
\end{eqnarray}
\begin{eqnarray}
\De_j(x):=(u^+_j-u^-_j)(x,p_j(x)), \ \ j=1,2,\ld\;.  \label{eq31}
\end{eqnarray}
Set
\begin{eqnarray}
u^-_i(x,t):=\int^t_{a_i(x)}\frac{ds}{\sqrt{h^2(x,s)+1}}\!+\!g^-_i(x),\; g^-_i(x)\!=\!\int^x_0\frac{h(s,a_i(s))\!+\!a'_i(s)}{\sqrt{h^2(s,a_i(s))\!+\!1}}ds,\;\text{in}\;\OO^-_i \hspace*{-2cm} \label{eq32}
\end{eqnarray}
\begin{eqnarray}
u_k(x,t)=\left\{\begin{array}{lll}
                  u^+_k(x,t), & \text{on} & \OO^+_k \\
                  u^-_k(x,t)+\De_k(x), & \text{on} & \OO^-_k
                \end{array}\right.,  \ \ k=1,2,\ld\;.  \label{eq33}
\end{eqnarray}
Next we define $u$ below $a_0(x)$.
\begin{eqnarray}
u^+_{-1}(x,t)=u^+_1(x,t) \ \ \text{on} \ \ \OO^+_{-1}=\{p_{-1}(x)\le t\le a_0(x)\},  \label{eq34}
\end{eqnarray}
with
\begin{eqnarray}
u^-_{-1}(x,t):=\int^t_{a_{-1}(x)}\frac{ds}{\sqrt{h^2(x,s)+1}}+g^-_{-1}(x), \;\text{on}\;\OO^-_{-1}=\{a_{-1}(x)\le t\le p_{-1}(x)\}\hspace*{-2cm}   \label{eq35}
\end{eqnarray}
where $g^-_{-1}(x)=\int\limits^x_0\dfrac{h(s,a_{-1}(s))+a'_{-1}(s)}{\sqrt{h^2(s,a_{-1}(s))+1}}ds$,
\begin{eqnarray}
\De_{-1}(x):=u^+_{-1}(x,p_{-1}(x))-u^-_{-1}(x,p_{-1}(x)),  \label{eq36}
\end{eqnarray}
\begin{eqnarray}
u_{-1}(x,t)=\left\{\begin{array}{lll}
                   u^+_{-1}(x,t), & \text{in} & \OO^+_{-1} \\
                   u^-_{-1}(x,t)+\De_{-1}(x), & \text{in} & \OO^-_{-1}
                   \end{array}\right.  .  \label{eq37}
\end{eqnarray}
And further down $i=2,3,\ld$ set
\begin{eqnarray}
\OO^+_{-i}=\{p_{-i}(x)\le t\le a_{-i+1}(x)\}, \ \ \OO^-_{-i}=\{a_{-i}(x)\le t\le p_{-i}(x)\},  \label{eq38}
\end{eqnarray}
\begin{eqnarray}
u^+_{-i}(x,t):=u^-_{-i+1}(x,t)+\De_{-i+1}(x), \ \ \text{on} \ \ \OO^+_{-i},  \label{eq39}
\end{eqnarray}
\begin{eqnarray}
\De_{-i}(x):=(u^+_{-i}-u^-_{-i})(x,p_{-i}(x)),  \label{eq40}
\end{eqnarray}
with
\begin{eqnarray}
u^-_{-i}(x,t):=\int^t_{a_{-i}(x)}\frac{ds}{\sqrt{h^2(x,s)+1}}+g^-_{-i}(x), \ \ \text{on} \ \ \OO^-_{-i}  \label{eq41}
\end{eqnarray}
where $g^-_{-i}(x)=\int\limits^x_0\dfrac{h(s,a_{-i}(s))+a'_{-i}(s)}{\sqrt{h^2(s,a_{-i}(s))+1}}ds$,
\begin{eqnarray}
u_{-k}(x,t)=\left\{\begin{array}{lll}
                     u^+_{-k}(x,t), & \text{in} & \OO^+_{-k} \\
                     u^-_{-k}(x,t)+\De_{-k}(x), & \text{in} & \OO^-_{-k}
                   \end{array}\right., \ \ k=2,3,\ld\;.  \label{eq42}
\end{eqnarray}
Finally set
\begin{eqnarray}
u(x,t)=u_k(x,t) \ \ \text{on} \ \ \OO^+_k\cup\OO^-_k, \ \ k\in Z\sm\{0\}.  \label{eq43}
\end{eqnarray}
With this definition we note that $u(x,t)$ is differentiable on $\R^2\sm S$, and
\begin{eqnarray}
\bigg(\frac{\partial u}{\partial t}\bigg)^2+\bigg(\frac{\partial u}{\partial x}\bigg)^2=1 \ \ \text{on} \ \ \R^2\sm S.  \label{eq44}
\end{eqnarray}
and thus we conclude as before that $h(x,t)\equiv$\,constant and $S=\emptyset$. The proof is complete.  \hfill$\square$
\end{Proof}

\end{document}